\documentclass[11pt,twoside]{article}
\usepackage{latexsym,amssymb,amsmath}
\usepackage{amsfonts}
\usepackage{color}
\definecolor{Gris}{cmyk}{0.1,0.1,0.1,.75}
\usepackage[colorlinks=true,citecolor=Gris,linkcolor=Gris,filecolor=Gris,urlcolor=Gris]{hyperref}
\RequirePackage{graphics,epsfig,psfrag}%,picins}
\def\abs#1{\left \vert #1 \right \vert}

\def\ZZ{{\bf Z}} %reelle Zahlen
\def\QQ{{\bf Q}} %reelle Zahlen
\def\cS{{\cal S}} %reelle Zahlen
\def\cO{{\cal O}} %reelle Zahlen
\def\Mod#1{\,(\hbox{\rm mod}\,#1)}

\def\TT{{\rm T}}

\def\Re{\hbox{\rm Re}\,}
\def\tr{\hbox{\rm tr}\,}

\def\phi{\varphi}
\def\eps{\varepsilon}

\def\undemi{{1\mskip-3mu /2}}
\def\pn{\medskip\par\noindent}
\def\bi{\vspace{-2pt}\begin{itemize}\itemsep -2pt plus 1pt minus 1pt}
\def\ei{\end{itemize}\vspace{-4pt}}
\def\bn{\vspace{-2pt}\begin{enumerate}\itemsep -2pt plus 1pt minus 1pt}
\def\en{\end{enumerate}\vspace{-4pt}}
\newcommand{\Pf}{{\em Proof}. }

\newcommand{\EPf}{\hbox{}\hfill$\Box$\vspace{.5cm}}
\def\[#1\]{\begin{eqnarray}#1\end{eqnarray}}
\def\$#1\${\begin{eqnarray*}#1\end{eqnarray*}}
\def\pent#1#2{\lfloor\frac{#1}{#2}\rfloor}

\def\comb#1#2{{}{\textstyle{{#1} \choose {#2}}}}

\def\Sum{\mathop{\sum}\limits}

\def\abs#1{\left \vert #1 \right \vert}

\def\frac#1#2{{\textstyle{{#1} \overwithdelims.. {#2}}}}
\def\Frac#1#2{{\displaystyle{{#1} \overwithdelims.. {#2}}}}
\catcode`\@=11
\def\system#1{\left\{\null\,\vcenter{\openup\jot\m@th
\ialign{
\strut\hfil$\displaystyle{##}$&
$\,\displaystyle{{}##}\,$\hfil&&
\strut\hfil$\,\displaystyle{##}$&
$\,\displaystyle{{}##}\,$\hfill
\hfil\crcr#1\crcr}}\right.}

\def\cmatrix#1{\left [
\null\,\vcenter{%\openup\jot\m@th
\ialign{
\hfil${##}\ $\hfil &
\hfil$\ {##}\ $\hfil&&
\hfil$\ {##}\ $\hfil&
\hfil$\ {##}$\hfil
\crcr#1\crcr}}\right ]}

\def\@begintheorem#1#2#3{\par\addvspace{8pt plus3pt minus2pt}%
              \noindent{\csname#1headfont\endcsname#1\ \ignorespaces#3 #2.}%
              \csname#1font\endcsname\hskip6pt\ignorespaces}
\def\@endtheorem{\par\addvspace{8pt plus3pt minus2pt}\@endparenv}

\catcode`\@=12

\newtheorem{theorem}{Theorem}[section]
\newtheorem{thm*}{Theorem}
\newtheorem{thm}[theorem]{Theorem}
\newtheorem{corollary}[theorem]{Corollary}
\newtheorem{lemma}[theorem]{Lemma}
\newtheorem{proposition}[theorem]{Proposition}

\newtheorem{remark}[theorem]{Remark}
\newtheorem{remarks}[theorem]{Remarks}
\newtheorem{example}[theorem]{Example}

\newtheorem{algorithm}[theorem]{Algorithm}
\newtheorem{conjecture}[theorem]{Conjecture}
\def\[#1\]{\begin{align}#1\end{align}}
\date{\today}
\title{On Alexander-Conway polynomials of two-bridge links}
\author{Pierre-Vincent Koseleff, Daniel Pecker}
\begin{document}
\maketitle
\begin{abstract}
We consider Conway polynomials of two-bridge links as Euler continuant
polynomials.
As a consequence, we obtain
new and elementary proofs of classical Murasugi's 1958 alternating theorem and Hartley's 1979 trapezoidal theorem.
We give a modulo 2 congruence for links, which implies the classical Murasugi's 1971 congruence for knots.
We also give sharp bounds for the coefficients of Euler continuants and deduce bounds for the
Alexander polynomials of two-bridge links.
These bounds improve and generalize those of Nakanishi-Suketa'96.
We  easily obtain some bounds for the roots of the Alexander
polynomials of two-bridge links. This is a partial answer to Hoste's conjecture on the roots of Alexander polynomials of alternating knots.
\pn
\noindent {\bf MSC2010}: 57M25, 11C08
\pn
{\bf Keywords}: Euler continuant polynomial, two-bridge link, Conway polynomial, Alexander polynomial
\end{abstract}
\begin{center}
\parbox{12cm}{\small
\tableofcontents }
\end{center}
\section{Introduction}
In this paper, we consider the Conway polynomial of a two-bridge link as an
Euler continuant polynomial.
We study the problem of determining whether a given polynomial is the Conway
polynomial of a two-bridge link (or knot), or equivalently, if it is a Euler continuant polynomial.
For small degrees, this problem can be solved by an exhaustive search
of possible two-bridge links.
Here, we give necessary conditions on the coefficients of the polynomial, which can be tested for high degree polynomials.
\pn
In section 2 we present Euler continuant polynomials and give some properties of their coefficients. We  show their relations with the
Fibonacci polynomials $f_k$ defined by:
$$
f_0=0, f_1=1, f_{n+2} (z)= z f_{n+1}(z) + f_n(z).
$$
In section 3, we recall the definitions of two-bridge links and we present the %Siebenmann's
description of the Conway polynomial of a two-bridge link as an extended Euler continuant polynomial.
%\pn
We obtain a characterization of modulo 2 two-bridged Conway polynomials.
\pn
{\bf Theorem \ref{fib2bis}.}
{\em
Let $\nabla (z) \in \ZZ[z]$ be the Conway polynomial of a rational  link (or knot).
There exists a Fibonacci polynomial $f_D(z)$ such that $ \nabla(z) \equiv f_D(z) \Mod 2$.}
\pn
We give a simple method (Algorithm \ref{degree}) that determines the integer $D$
such that $ \nabla(z) \equiv f_D(z) \Mod 2$.
This is used to test when $\nabla (z) \equiv 1 \Mod 2,$
which is a necessary condition to be a two-bridge Lissajous knot.
\pn
These results are applied in section 4 to the Conway polynomials of two-bridge links
denoted
$$ \nabla_m (z) = \sum_{k=0}^{\pent m2 } c_{m-2k} z^{m-2k}.$$
%\pn
{\bf Theorem \ref{nak1bis}.}
{\em
For $ k\geq 0$,
$$\abs {c_{m-2k}} \le   \comb{ m-k}{k }   \abs {c_m}.$$
If equality holds for some  positive integer $k< \pent m2$,
then it holds for all integers. In this case,  the link is  isotopic to a link
 of Conway form $ C ( 2, -2, 2,  \ldots, (-1)^{m+1} 2 )$
 or $ \allowbreak C(2,2, \ldots, 2) $, up to mirror symmetry.}
\pn
When $ |c_m| \ne 1,$ we have the following sharper bounds:
\pn
{\bf Theorem \ref{nak3bis}.}
{\em Let $ g \ge 1$ be the greatest prime divisor of $c_m, $ and $m\ge 2k \ge 2$. Then
$$
\abs { c_{m-2k}} \le \Bigl( \comb{ m-k-1}{k}  + \Frac 1g \bigl( \comb{ m-k-1}{k-1} -1 \bigr)
\Bigr) \abs {c_m} +1.
$$
Equality holds for links of Conway forms
$C (2g, 2,2, \ldots, 2 )$  and $ C( 2g, -2,2, \ldots, (-1)^{m +1}\, 2)$.}
\pn
In section 5, we apply our results to the Alexander polynomials.
Our modulo 2 congruence of Theorem \ref{fib2bis} provides  a simple proof of a congruence of Murasugi \cite{Mu2} for periodic knots (two-bridge knots have period two).
Moreover, we deduce a  congruence for the Hosokawa polynomials of two-bridge links (Corollary \ref{hosokawa}).
\pn
Then, we obtain a simple proof of both the Murasugi alternating theorem \cite{Mu,Mu1}, and the
Hartley trapezoidal theorem \cite{Ha} (see also \cite{Kan1}) using the trapezoidal property:
\pn
{\bf Theorem \ref{alphabis}.}
{\em Let $K $
be a two-bridge  link (or knot).
Let
$$
\nabla_K=
c_m
 \Bigl(\sum_{i=0}^{\pent m2} (-1)^i \alpha_{i} f_{m-2i+1}\Bigr),  \ {}   \
\alpha_0=1
 $$
be its Conway polynomial written in the Fibonacci basis. Then we have
\bn
\item
$\alpha_{j}\geq 0, \ j = 0,\ldots,\pent m2$.
\item If $\alpha_{i}=0$ for some $i>0$ then
$\alpha_{j}=0$ for $j\geq i$.
\en
}
\pn
We conclude this section with bounds for the coefficients of the
Alexander coefficients. These  bounds improve those of Nakanishi and Suketa for the
Alexander polynomials of two-bridge knots (see \cite[theorems~2 and 3]{NS}).
Moreover, they are sharp and  hold for any $k$.
\pn
We prove that the conditions on Conway coefficients are sharper than the conditions
on the Alexander coefficients deduced from them.
\pn
In section 6, we conclude our paper
with the following convexity conjecture:
\pn
{\bf Conjecture \ref{con-con}}.
{\em Let $\Delta(t) = a_0 -a_1 ( t+ t^{-1}) + a_2 ( t^2 + t^{-2} ) - \cdots +(-1)^n a_n ( t^n + t^{-n} )$ be
the Alexander polynomial of a two-bridge  knot. Then there exists an integer $k \le n$ such that
$ (a_0, \ldots, a_k) $ is convex and $(a_{k}, \ldots, a_n)$ is concave.}
\pn
We have tested this conjecture for all two-bridge knots with 20 crossings or fewer.
\pn
We also deduce some bounds for the roots of Alexander polynomials of two-bridge links (or knots) from the properties of Euler continuant polynomials. This gives some partial answer to the Hoste conjecture \ref{hoste}.
\section{Extended Euler continuant polynomial}\label{cf}
We define the extended Euler continuant polynomial $D_m(b_1,\ldots,b_m)(z)$ as the determinant of the tridiagonal matrix
\[
\left (
\begin{array}{ccccc}
b_1 z & -1 & 0 & \ldots & 0\\
1 & b_2 z & -1 & \ddots &\vdots \\
0 & & \ddots&\ddots & 0 \\
\vdots & \ddots & \ddots & & -1 \\
0 & \ldots & 0 & 1 & b_m z
\end{array}
\right ). \label{tridiagonal}
\]
The polynomials $D_i$ satisfy the recurrence relation
\[
D_{-1} = 0,\, D_0 = 1, \, D_k = b_k z D_{k-1} + D_{k-2}.
\]
When $z=1,$ this is the  classical Euler continuant polynomial (see \cite{Kn}).
\pn
When all the $b_i$ are equal to 1, we obtain the Fibonacci polynomials defined by
\[ f_0=0, f_1=1, f_{n+2} (z)= z f_{n+1}(z) + f_n(z),   \ n \in \mathbf{Z} . \label{fib}\]
Let us recall some basic facts about Fibonacci polynomials.
\begin{lemma}\label{fibex}
For $m \ge 0$:
$$
f_{m+1}(z) = \sum_{k=0}^{ \pent m2}
\comb{ m-k}{k} z ^{m-2k}.
$$
\end{lemma}
\Pf By induction on $m$. The result is clear for $m=1$ and for
$m=2$.
Let us suppose the result true for $m-1$
and $m$.
By induction, the coefficient of $ z^{m-2k}$  is
 $\comb{m-1-k}{k}$ in $z f_m(z),$ and $ \comb{m-1-k}{k-1} $ in $f_{m-1} (z)$.
Consequently, the coefficient of $z^{m-2k}$ in $f_{m+1}(z)$ is\\
\hbox{} \hfill
$
\comb{ m-1-k}{k} + \comb{ m-1-k}{k-1} = \comb{m-k}{k}.
$\hfill
\EPf
%\pn
%%%%%%%%%%%%%%%%%%%%%%%%%%%%%%%%%%%%%%%%%%%%%%%%%%%%%%%%%%%
\begin{remark}
This means that the Fibonacci polynomials can be read on the  diagonals of  Pascal's triangle.
When $z=1,$ we recover the classical Lucas identity
$$
F_m = \sum_{k=0}^{ \pent m2} \comb{ m-k}{k},
$$
where $F_m$ are  the Fibonacci numbers ($F_0=0, \  F_1=1, \ldots, F_{n+1}=F_n + F_{n-1}$).
\end{remark}
%\section{Inequalities for the coefficients of Euler continuant polynomial}
\pn
We shall need the following  explicit notation for Euler continuant polynomials:
\[
D_m (z)= \sum _{k=0}^{ \pent m2} c_{m-2k} (b_1, \ldots, b_m) z^{m-2k}.
\]
We obtain some properties of $ c_{m-2k} (b_1, \ldots, b_m),$
considered as a polynomial in the $m$ variables $b_1, \ldots,  b_m.$
\begin{proposition}\label{euler}
Let ${\cal M}$ be the set of all monomials
$\Frac {b_1 \cdots b_m}{ b_{i_1}b_{i_1+1} \cdots b_{i_k} b_{i_k +1}},$
where $k \ne 0$ and $i_h +1 < i_{h+1}$. Let ${\cal M}_j$ be the subset of all monomials
of ${\cal M}$ that are relatively prime to $b_j.$
Then we have
\bn
\item  The polynomial
$c_{m-2k} (b_1, \ldots, b_m) $ is the sum of all monomials of ${\cal M}$.
\item  The set ${\cal M}$ has
 $\comb{ m-k}{k}$ elements.
\item The monomials of ${\cal M}$ do not have a common divisor except 1.
\item The number of elements of ${\cal M}_j$ is at least
 $ \comb{ m-1-k}{k-1}$.
\item If $m \ge 4$, then the monomials of ${\cal M}_j$ do not have a common divisor except 1.
\en
\end{proposition}
\Pf
\bn
\item
This is a classical property of the Euler continuant (see \cite{Kn})
\item This number is  $c_{m-2k} (1,1, \ldots, 1),$
which is a coefficient of the Fibonacci polynomial
$$f_{m+1} (z)=
\sum_ {k=0}^{ \pent m2}c_{m-2k} (1,1, \ldots, 1)z^{m-2k}
= \sum_ {k=0}^{\pent m2} \comb{ m-k}{k} z^{m-2k}.$$
\item For every integer $i \le m,$ there is an element of ${\cal M}$ which is not divisible
by  $b_i.$
Hence the GCD of the elements of ${\cal M}$ is 1.
\item
Let $1\leq j \leq m$ and $\mathbf{b} = (1,\ldots,1,0,1,\ldots,1)$ where $b_j=0$, and $b_k=1$ for $k\ne j$.
Let us define the polynomials $g_n, $ for $n \le m$ by
$g_n(z)= D_n (\mathbf{b}) (z)$.
The number of  elements of ${\cal M}_j$  is the coefficient
$ c_{m-2k}(\mathbf{b})$ of $g_m(z)$.
%Let us compute the first polynomials $g_n(z)$.
\pn
If $j=1$, then we have $g_1=0$, $g_2=1.$
Then, an easy induction shows that  $ g_n= z g_{n-1} + g_{n-2}$ is the Fibonacci polynomial
 $g_n=f_{n-1}. $
\pn
If $j>1$, then we have
$$
g_1 = f_2, \  \ldots, \  g_{j-1}= f_{j}, \   g_j = f_{j-1}, \quad
{\rm and} \quad  g_{n+1}= z g_n + g_{n-1} \quad {\rm if}  \quad  n \geq j.$$
Let us  write $ p(z)\succeq q(z)$ when each coefficient of $p$ is
greater than or equal to the corresponding coefficient of $q$.
We have $f_{k+2} \succeq f_k$, and therefore
$ g_{j+1}= z f_{j-1} + f_j \succeq z f_{j-1} + f_{j-2}= f_j .$
Then a simple induction shows that $g_m \succeq f_{m-1}$, and consequently
$ c_{m-2k}( {\bf b}) \ge \comb {m-1-k}{k-1}.$
\item
Since $ m \ge 4,$  for every  $i \ne j,$ there is
a monomial which is not divisible by  $b_i$. Consequently, the GCD of the elements of
$ {\cal M}_j$ is 1.
\EPf
\en
%%%%%%%%%%%%%%%%%%%%%%%%%%%%%%%%%%%%%%%%%%%%%%%
\section{Conway polynomials of two-bridge links}
\begin{figure}[th]
\begin{center}
{\scalebox{.8}{\includegraphics{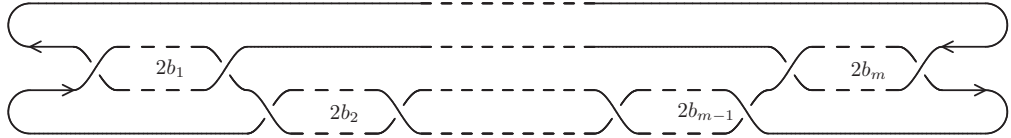}}}
\end{center}
\caption{Oriented two-bridge links ($m$ odd)}
\label{conways2}
\end{figure}

A two-bridge knot (or link) admits a diagram in Conway's normal form.
This form, denoted by
$C(a_1, a_2, \ldots, a_n)$  where $a_i$ are integers, is explained by the following
picture (see \cite{Co,Mu}).
\psfrag{a}{\small $a_1$}\psfrag{b}{\small $a_2$}%
\psfrag{c}{\small $a_{n-1}$}\psfrag{d}{\small $a_{n}$}%
\begin{figure}[h!]
\begin{center}
{\scalebox{.8}{\includegraphics{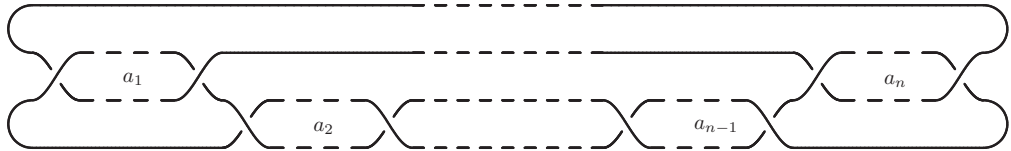}}}
%\\[30pt]
%{\scalebox{.8}{\includegraphics{kr1bb.eps}}}
\end{center}
\caption{Conway's normal forms}
\label{conways1}
\end{figure}
\pn
The number of twists is denoted by the integer
$\abs{a_i}$, and the sign of $a_i$ is defined
as follows: if $i$ is odd, then the right twist is positive,
if $i$ is even, then the right twist is negative.
On  Fig. \ref{conways1} the $a_i$ are positive (the $a_1$ first twists are right twists).
\pn
The two-bridge links are classified by their Schubert fractions (see \cite{Sc})
$$
\Frac {\alpha}{\beta} =
a_1 + \Frac{1} {a_2 + \Frac {1} {\cdots +\Frac 1{a_n}}}=
[ a_1, \ldots, a_n], \quad \alpha >0.
$$
We shall denote  $S \bigl( \Frac {\alpha}{\beta} \bigr)$  a two-bridge link with
Schubert fraction $ \Frac {\alpha}{\beta}.$
The two-bridge links
$ S (\Frac {\alpha} {\beta} )$ and $ S( \Frac {\alpha ' }{\beta '} )$ are equivalent
if and only if $ \alpha = \alpha' $ and $ \beta' \equiv \beta ^{\pm 1} ( {\rm mod}  \  \alpha).$
The integer $ \alpha$ is odd for a knot, and even for a two-component link.
\pn
When $\alpha\beta$ is even, one shows (see \cite[p.~26]{Kaw},
\cite{Ko,KM}) that there is a unique continued fraction expansion
$\Frac{\alpha}{\beta} = [2b_1, 2b_2, \ldots, 2b_n], \ b_i \in \ZZ - \{0\}$.
It means that
any oriented two-bridge link can be put in the form shown in Figure \ref{conways2}.
It will be denoted by $C(2b_1, 2b_2, \ldots,2\, b_m)$,
including the indicated orientation.
This is a two-component link if and only if $m$ is odd.
\psfrag{a}{\small $2b_1$}\psfrag{b}{\small $2b_2$}%
\psfrag{c}{\small $2b_{m-1}$}\psfrag{d}{\small $2b_{m}$}%
\pn
The Conway polynomial $\nabla_K(z) \in \ZZ[z]$
is a polynomial invariant of the oriented link $K$ (see \cite{Cr}). When $K$ is a two-bridge link
its Conway polynomial $\nabla _m$ is given by the following
method (see \cite{Si} and \cite[Th. 8.7.4]{Cr}):
\begin{thm}[\cite{Si,Cr}]\label{th:siebenmann}
Let us consider the oriented two-bridge link
$$C( 2b_1, -2b_2, \ldots, (-1)^{m-1} 2 b_m).$$
Its Conway polynomial $\nabla_m(z)$ is the Euler continuant polynomial
$D_m(b_1,\ldots,b_m)(z)$.
\end{thm}
\begin{example}[The torus links]
{\it The Conway polynomial of the torus link $\TT(2,m)$ is the Fibonacci polynomial
$f_{m}(z)$ (see  \cite{Ka,KP4})}.
\end{example}
Consequently, the following result gives in fact a characterization of
modulo 2 Conway polynomials of two-bridge links.
\begin{theorem}\label{fib2bis}
Let $\nabla (z) \in \ZZ[z]$ be the Conway polynomial of a rational  link (or knot).
There exists a Fibonacci polynomial $f_d(z)$ such that $ \nabla(z) \equiv f_d(z) \Mod 2$.
\end{theorem}
\Pf
Let us  write
$ (a,b) \equiv (c,d) \Mod 2 $ when $ a \equiv c \Mod 2$ and $b \equiv d \Mod 2$.
We will show by induction on $m$ that there exist integers $d$ and $e = \pm 1$ such that
$(D_{m-1},D_{m}) \equiv (f_{d-e},f_{d}) \Mod 2$.
\pn
The result is true for $m=0$ as $(\nabla_{-1},\nabla_0) = (0,1) = (f_0,f_1)$, that is $d=e=1$.
\pn
Suppose that  $(\nabla_{m-1},\nabla_{m}) \equiv (f_{d-e},f_{d}) \Mod 2$, with $e=\pm 1$
for some $m\geq 0$.
Then we have $\nabla_{m+1} = b_{m+1} z \nabla_{m} + \nabla_{m-1}$.
\pn
If $b_{m+1}\equiv 0 \Mod 2$ then
$
\nabla_{m+1} \equiv  \nabla_{m-1} \equiv  \nabla_{d-e} \Mod 2
$
and $(\nabla_{m},\nabla_{m+1}) \equiv (f_{d},f_{d-e})$.
\pn
If $b_{m+1}\equiv 1 \Mod 2$ then
$
\nabla_{m+1} \equiv  z f_{d} + f_{d-e} \equiv f_{d+e} \Mod 2.$
Consequently $(\nabla_{m},\nabla_{m+1}) \equiv (f_{d},f_{d+e})$.
\EPf
\pn
We thus deduce a fast algorithm
for the determination of the integer $d$ such that $\nabla_m \equiv f_d \Mod 2$,
see also \cite{Bu}.
%%%%%%%%%%%%%%%%%%%%%%
\begin{algorithm}\label{degree}
%Let $(b_1, \ldots, b_m)$.
Let us define the sequences of integers $e_i$ and $d_i$, $ i=0,\ldots, m,$ by
$$e_0=1, \  d_0 = 1, \ e_{i+1} = -(-1)^{b_{i+1}}e_i, \
d_{i+1} = d_i + e_{i+1}.$$
Then we have  $\nabla_m(z) \equiv   f_d (z) \Mod 2$ where $d= \abs {d_m}$.
\end{algorithm}
%%%%%%%%%%%%%%%%%%%%%%%%%%%%%%%%%%%%%%%%%%%%%%%%%%%%%%%%%%
\begin{remark} Let us consider the two-bridge link $K=C(2b_1, -2b_2, \ldots, (-1)^{m-1} 2 b_m)$.
From \cite{St07}, the crossing number $N$ of $K$ is $2 \sum_{i=1}^{m} \abs{b_i} - \#\{i, b_i b_{i+1} <0\} \geq m+1$.
We deduce that one computes $d$ such that $\nabla_K \equiv f_d \Mod 2$ in $\cO(N)$ steps.

The torus knot ${\TT}(2,m)$ is the two-bride knot $S(m)$ of crossing number $m$.
The rational number $\frac{m}{m-1}$ has the continued fraction expansion of length $m-1$: $[2,-2, \ldots, (-1)^{m} 2]$. That shows that the inequality $m \leq N-1$ is sharp.
\end{remark}
%%%%%%%%%%%%%%%%%%%%%%%%%%%%%%%%%%%%%%%%%%%%%%%%%%%%%%%%%%%%%%
Jones, Przytycki, and Lamm proved that the Conway polynomial of a
two-bridge Lissajous knot satisfies the congruence
$ \nabla (z) \equiv 1 \Mod2,$ that is $d=0$ (see \cite{JP,La}).
Using Algorithm \ref{degree}
we deduce the number of two-bridge knots with a Conway polynomial
congruent to 1 modulo 2 (see Table \ref{lisst} and compare \cite{BDHZ}).
\begin{table}[th] %\label{derniere}
\begin{center}
{\small
\setlength\tabcolsep{2pt}
\newlength\tot
\tot=\textwidth \advance\tot by -8\tabcolsep % 2 * nombre de colonnes
\begin{tabular}{||p{.20\tot}|p{.040\tot}|p{.040\tot}|p{.050\tot}|p{.050\tot}|p{.050\tot}|p{.070\tot}|p{.075\tot}|p{.075\tot}|p{.075\tot}|p{.075\tot}|p{.075\tot}|p{.075\tot}||}
\hline
\hline
Crossing Number &\hfill 3&\hfill 4&\hfill 5&\hfill 6&\hfill 7&\hfill 8&\hfill 9&\hfill 10&\hfill 11&\hfill 12\\
\hline
Two-bridge &\hfill 1&\hfill 1&\hfill 2&\hfill 3&\hfill 7&\hfill 12&\hfill 24&\hfill 45&\hfill 91&\hfill 176\\
\hline
$\nabla(t)\equiv 1$ &\hfill 0&\hfill 0&\hfill 1&\hfill 1&\hfill 2&\hfill  4&\hfill  8&\hfill 13&\hfill 26&\hfill  51\\
\hline
\end{tabular}
\vspace{20pt}\\
\begin{tabular}{||p{.20\tot}|p{.040\tot}|p{.040\tot}|p{.050\tot}|p{.050\tot}|p{.050\tot}|p{.070\tot}|p{.075\tot}|p{.075\tot}|p{.075\tot}|p{.075\tot}|p{.075\tot}|p{.075\tot}||}
\hline
Crossing Number&\hfill 13&\hfill 14&\hfill 15&\hfill 16&\hfill 17&\hfill 18&\hfill 19&\hfill 20&\hfill 21&\hfill 22\\
\hline
Two-bridge &\hfill 352&\hfill 693&\hfill 1387&\hfill 2752&\hfill 5504&\hfill 10965&\hfill 21931&\hfill 43776&\hfill 87552&\hfill 174933\\
\hline
$\nabla(t)\equiv 1$ &\hfill  97&\hfill 185&\hfill  365&\hfill  705&\hfill 1369&\hfill  2675&\hfill  5233&\hfill 10211&\hfill 20011&\hfill  39221\\
\hline
\hline
%\end{tabular}}
\end{tabular}}
\end{center}
\vspace{-5pt}
\caption{The number of two-bridge knots, and two-bridge
 knots with Conway polynomial congruent to 1 modulo 2.}\label{lisst}
\end{table}
%%%%%%%%%%%%%%%%%%%%%%%%%%%%%%%%%%%%%
%%%%%%%%%%%%%%%%%%%%%%%%%%%%%%%%%%%%
\section{Inequalities for Conway Polynomials}
We shall write the Conway polynomial of a two-bridge link
$$ \nabla_m (z) = \sum_{k=0}^{\pent m2 } c_{m-2k} z^{m-2k}.$$
%Theorems \ref{nak1} and \ref{nak3} on Euler continuant polynomials give:
\begin{thm}\label{nak1bis}
For $ k\geq 0$,
$$\abs {c_{m-2k}} \le   \comb{ m-k}{k }   \abs {c_m}.$$
If equality holds for some
%positive
integer $k< \pent m2$,
then it holds for all integers. In this case,  the link is  isotopic
to the torus link $T(2,m) $
or to the link $ C(2,2, \ldots, 2) $, up to mirror symmetry.
\end{thm}
\Pf
By Proposition \ref{euler},  the number of monomials of
$c_{m-2k} (b_1, \ldots, b_m)$  is
$ \comb{ m-k}{k}$. The result follows since no monomial
 is greater than $\abs{c_m}= \abs{ b_1\cdots b_m}$.
\pn
If equality holds for some positive integer $k< \pent m2$,  then for all $i, \, j, \  $
$ b_i b_{i+1}= b_j b_{j+1} = \pm 1,$ which implies the result.
\EPf
\begin{example}
The knot $10_{145}$ has Conway polynomial $P=1+5z^2+z^4$. We have $P
\equiv f_5 \Mod 2,$ but $P$ does not satisfy the condition $\abs{c_2}
\leq 3$, and then $10_{145}$ is not a two-bridge knot.
\pn
The knot $11n109$ has Conway polynomial $1+6z^2+z^4 - z^6 .$
%= f_7 - f_3+ f_1.
It satisfies the bounds of Theorem \ref{nak1bis}: $\abs{c_2} \leq 6, \
\abs{c_4} \leq 5$, but not the equality condition:  $c_2=6$ whereas $c_4 \ne 5$.
Consequently, $11n109$ is not a two-bridge knot.
\end{example}
%%%%%%%%%%%%%%%%%%%%%%%%%%%%%%%%%%%%%%%%%%%%
%To prove the refined inequalities of Theorem \ref{nak3}, we shall use
We shall use the following lemma, which generalizes the inequality
$ a+b \le ab+1,$  valid for positive integers (see also \cite{NS}).
\begin{lemma}\label{fact1}
Let $p_i, \, i \in \cS $  be relatively prime  divisors of
$p = x_1 x_2 \cdots x_m$ in $\QQ[x_1,\ldots, x_m]$.\\
Let  $ {\bf b}= (b_1, \ldots, b_m) $ be a m-tuple of positive integers.
Then
\[
\sum _{i \in \cS} p_i ( {\bf b})  \le
\Bigl( {\rm card} (\cS)-1 \Bigr) {p( {\bf b})} +1.\label{fact2}
\]
\end{lemma}
\Pf
We do not suppose the $p_i$ distinct.
Let us prove the result by induction on $k={\rm card}(\cS)$.
\medskip
The result is clear if $k=1,$ we have ${p_1}=\pm 1,$ and
the inequality is $\pm 1 \le 1$.
\pn
If all the $p_i=1,$ the result is clear.
Otherwise, let $x_h$ be a divisor of some $p_i$.
\pn
Let $\cS_1= \{ i\in \cS : \   x_h |p_i \}, $ and $\cS_2= \cS-\cS_1$.
We have $k=k_1+k_2$,
%\ k_1< k, \  k_2 < k$,
where $k_j= {\rm card} (\cS_j)$.
Let $q_j= {\rm GCD}\{ p_i, i \in \cS_j \}$, then $q_1$ and $q_2$ are coprime,
and $q_1q_2$ is a divisor of $p$.
\pn
By induction we obtain for $j= 1, 2$:
$$ \sum _{i \in \cS_j} p_i( {\bf b}) \le {q_j( {\bf b})}
\Bigl( ( k_j-1) \Frac { {p( {\bf b} )}} {{q_j ( {\bf b})} } +1 \Bigr)
= (k_j-1)p ({\bf b}) + q_j( {\bf b }) .$$
Adding these two inequalities we get
\begin{align*}
\sum _{ i \in \cS} p_i ( {\bf b}) &\le
( k_1+k_2-1) {p( {\bf b})} + {q_1 ( {\bf b})} + {q_2 ( {\bf b})} - {p({\bf b})}\\
&\le
(k_1+k_2-1) {p({\bf b})} + {q_1({\bf b})} {q_2( {\bf b})}  - {p({\bf b})}+1,
\end{align*}
which proves the result, since ${q_1({\bf b})}{ q_2( {\bf b})} \le {p({\bf b})}$.
\EPf
%%%%%%%%%%%%%%%%%%%%%%%%%%%%%%%%%%%%%%%%%%%%%%%%%%%%%%%%%%
\pn
With this lemma we can prove:
\begin{thm}\label{nak3bis}
Let $ g \ge 1$ be the greatest prime divisor of $c_m$, and let $k \ne 0.$ Then
$$
\abs { c_{m-2k}} \le \Bigl( \comb{ m-k-1}{k}
 + \frac 1g \bigl( \comb{ m-k-1}{k-1} -1 \bigr)
\Bigr) \abs {c_m} +1.
$$
Equality holds for %for links of Conway form
$(b_1, \ldots, b_m) = (g,1, \ldots, 1)$ and
$(b_1, \ldots, b_m) = (g,-1, \ldots, (-1)^m)$.%2g,-2, \ldots, (-1)^{m-1} 2)$ and $C(2g,2, \ldots,  2)$.
\end{thm}
\Pf%[\bf Proof of Theorem \ref{nak3}. ]
If $k=1$,
there are $m-1$ monomials in the polynomial
$c_{m-2}(b_1,\ldots,b_m)$, by Proposition \ref{euler}.
Then, using Lemma \ref{fact1} and the notation
$ \abs {\bf b} =(\abs {b_1}, \ldots, \abs {b_m} )$, we get
$$
\abs{c_{m-2}} = \abs {c_{m-2}({\bf b})} \le c_{m-2}( \abs {\bf b} )
 \leq (m-2) c_{m}(\abs{\bf b}) + 1 = (m-2) \abs{c_m} + 1.
$$
Now, suppose $ k \ge 2$. Let $g$ be the greatest prime divisor of
the integer $c_m= b_1 \cdots b_m,$ and suppose that $g \mid b_j$.
Let $N$ be the number of monomials of $c_{m-2k}( b_1, \ldots, b_m)$
that are prime to the monomial $b_j$. By Proposition \ref{euler},
these monomials are relatively prime, and $N \ge \comb{
m-1-k}{k-1}$. Using Lemma \ref{fact1} we obtain:
$ \sum_{ p_i \in
{\cal M}_j} p_i ({\bf b}) \le (N-1) \Frac{\abs{ c_m} }{\abs{b_j}}
+1$
and then
\begin{align*}
\abs {c_{m-2k}} = \abs {\textstyle{\sum _{ p_i \in {\cal M}}} p_i
({\bf b})}
&\le \Bigl( \Frac{N-1}{g}+ (\comb{m-k}{k}-N ) \Bigr)\abs {c_m}+1\\
&=
\Bigl( \comb{m-k}{k} -N (1 - \frac 1g) - \frac 1g\Bigr)\abs {c_m}+1\\
&\le
\Bigl( \comb{ m-k}{k} - \comb{ m-1-k}{k-1} (1-\frac 1g) - \frac 1g  \Bigr)\abs {c_m}+1\\
&=\Bigl( \comb{ m-1-k}{k} + \frac 1g  (\comb{ m-1-k}{k-1}-1 ) \Bigr)
\abs {c_m} +1.
\end{align*}
%For links  of Conway form $C(2g,-2, \ldots, (-1)^{m+1} 2)$, we have
For $\mathbf{b} = (g, 1, \ldots, 1)$ we obtain
$N=\comb{ m-1-k}{k-1}$, $c_m=g,$ and
$c_{m-2k}= g \comb {m-1-k}{k} + \comb {m-1-k}{k-1},$ and equality holds throughout.
\pn
For $\mathbf{b} = (g, -1, 1, \ldots, (-1)^{m})$
we get
$c_{m-2k}= (-1)^{{\pent m2} + k} \Bigl( g \comb {m-1-k}{k} + \comb {m-1-k}{k-1} \Bigr)$.
%which proves the result.
\EPf
\begin{example}
The knot $13n3010$ has Conway polynomial
$\nabla = 1+10\,{z}^{2}+4\,{z}^{4}-2\,{z}^{6}$.
It satisfies all conditions of Theorems \ref{nak1bis} and \ref{fib2bis} but not those of Theorem \ref{nak3bis}.
\end{example}
\pn
Now, we will express the Conway polynomials in terms of Fibonacci polynomials, and show that their coefficients are alternating.
\begin{thm}\label{alphabis}%
Let $K $ be a two-bridge  link (or knot).
Let
$$
\nabla_K= c_m
 \Bigl(\sum_{i=0}^{\pent m2} (-1)^i \alpha_{i} f_{m-2i+1}\Bigr),  \ {}   \
\alpha_0=1
 $$
be its Conway polynomial expressed in the Fibonacci basis. Then we have
\bn
\item
$\alpha_{j}\geq 0, \ j = 0,\ldots,\pent m2$.
\item If $\alpha_{i}=0$ for some $i>0$ then
$\alpha_{j}=0$ for $j\geq i$.
\en
\end{thm}
\Pf
We have $\nabla_0 = f_1$, $\nabla_1 = b_1 f_2$, $\nabla_2 = b_1b_2 \Bigl(f_3 - (1-\frac{1}{b_1b_2}) f_1\Bigr)$.

Let us show by induction that if
$$
\nabla_m = b_1 \cdots b_m \Bigl(\Sum_{i=0}^{\pent m2} (-1)^i \alpha_{i} f_{m+1-2i}\Bigr), \,
\nabla_{m-1} = b_1 \cdots b_{m-1} \Bigl(\Sum_{i=0}^{\pent{m-1}2} (-1)^i \beta_{i} f_{m-2i}\Bigr)$$
then $\alpha_{j} \geq  \beta_{j}\geq 0$, and if $\alpha_i=0$ for some $i,$ then
$\alpha_j=0$ for $j \ge i$.
\pn
The result is true for $m=2$ from the expressions of $\nabla_1$ and $\nabla_2$.
Using $zf_{m+1-2i} =  f_{m+2-2i}-f_{m-2i}$ and $\nabla_{m+1} = b_{m+1}z \nabla_m + \nabla_{m-1}$,
we deduce that
$$\nabla_{m+1} = b_1 \cdots b_{m+1} \Bigl(\Sum_{i=0}^{\pent {m+1}2} (-1)^i \gamma_{i} f_{m+2-2i}\Bigr),$$
where
$\gamma_{0} = 1$ and
\[
\gamma_{i}
&= \alpha_{i} + (\alpha_{i-1} - \beta_{i-1}) + (1- \frac{1}{b_m b_{m+1}}) \beta_{i-1},
\, i = 1,\ldots, \pent{m+1}2. \label{abc}
\]
As $\abs{b_m b_{m+1}}\geq 1$, we deduce by induction that $\gamma_{i} \geq \alpha_{i} \geq 0$.

Furthermore, if $\gamma_i=0,$ then by Formula (\ref{abc}) $ \alpha_i=0,$ and then, by induction,
$\alpha_j = \beta_j=0$ for $j\ge i$.   Finally, by Formula (\ref{abc}), we get $\gamma_j=0$
for $j \ge i$.
\EPf
%%%%%%%%%%%%%%%%%%%%%%%%%%%%%%%%%%%%%%%%%%%%%%%%%%%%%%%%%%%%%%%%%%%%%%%%%%%%%%%%%
\section{Applications to the Alexander polynomial}
In this section, we will see that our necessary conditions on the
Euler continuant polynomials imply analogous necessary conditions on both
Conway coefficients and Alexander coefficients of two-bridge knots and links. These
conditions are improvements of the classical results.
\pn
The Conway and the Alexander polynomials of a knot $K$ will be denoted by
$$\nabla_K (z) = 1 +{\tilde c}_1 z^2 + \cdots + {\tilde c}_n z^{2n}$$
and
$$ \Delta_K(t) = a_0 - a_1 (t+t^{-1}) + \cdots + (-1)^n a_n ( t^n+ t^{-n} ).$$
The Alexander polynomial $ \Delta_K (t)$ is deduced from the Conway polynomial:
$$ \Delta_K (t) = \nabla_K \Bigl( t^{\undemi} - t^{ - \undemi} \Bigr). $$
It is often normalized so that  $a_n$ is positive.
Thanks to this formula, it is not difficult to deduce
the Alexander polynomial from the Conway polynomial.
If we use the Fibonacci basis, it is even easier to deduce the
Conway polynomial of a knot from its Alexander polynomial.
\begin{lemma}\label{undemi}
If $ \ z =t^{\undemi} - t^{ -\undemi}, \ $  and $ \  n \in \mathbf{Z} \ $ is an integer,
then we have the identity
$$ f_{n+1} (z) + f_{n-1} (z) = (t^{\undemi})^n + ( -t^{ - \undemi} )^n,$$
where $f_k (z)$ are the  Fibonacci polynomials.
\end{lemma}
\Pf
Let
$ A= \cmatrix { z & 1 \cr 1 & 0 } $ be the  (polynomial) Fibonacci matrix.
If   $z =t^{\undemi} - t^{-\undemi}$, the eigenvalues of $A$
are  $ \  t^{\undemi}$ and $ \  -t^{ -\undemi}$, and consequently
$ \tr A^n =(t^{\undemi})^n + ( -t^{ -\undemi} )^n$.
On the other hand,
we have $ A^n= \cmatrix { f_{n+1}(z) & f_n(z) \cr f_n(z) & f_{n-1}(z) }$, and then
$ \tr A^n = f_{n+1} (z) + f_{n-1} (z)$.
\EPf
%%%%%%%%%%%%%%%%%%%%%%%%%%%%%%%%%%%%%%%%%%ש
\pn
From Lemma \ref{undemi}, we immediately deduce:
\begin{proposition} \label{chvar}
Let the Laurent polynomial $P(t)$ be defined by
$$P(t) = a_0 -a_1 (t+ t^{-1}) + a_2 ( t^2 + t^{-2} ) - \cdots +(-1)^n a_n ( t^n + t^{-n} ).$$
We have
$$P(t)= \sum_{k=0}^n (-1)^k (a_k-a_{k+1}) f_{2k+1}(z),$$
where $z =t^{\undemi} - t^{-\undemi}, $ and $a_{n+1}=0$.
\end{proposition}
We deduce a useful formula
(by substituting $a_0 = \ldots = a_n
=1$).
\[
f_{2n+1} \bigl( t^{\undemi} - t^{- \undemi} \bigr)=
(t^n + t ^{-n}) - (t^{n-1} + t ^{1-n}) + \cdots + (-1)^n.
\label{afib}
\]
Then, we deduce a simple proof of an elegant criterion due to
Murasugi (\cite{Mu2,Bu})
\begin{corollary}[Murasugi (1971)] {\label{al1}}
Let
$\Delta (t) = a_0 -a_1 ( t+ t^{-1}) + a_2 ( t^2 + t^{-2} ) - \cdots +(-1)^n a_n ( t^n + t^{-n} )$ be the Alexander polynomial of a two-bridge knot.
There exists an integer $k \le n$ such that
$ a_0, a_1, \ldots, a_k$ are odd, and $a_{k+1}, \ldots, a_n $ are even.
\end{corollary}
\Pf
If $K$ is a two-bridge knot, its Conway polynomial is a modulo 2 Fibonacci polynomial $f_{2k+1}$
by theorem \ref{fib2bis}. By Proposition \ref{chvar} we have
$
 f_{2k+1} \bigl( t^{\undemi} - t^{- \undemi} \bigr)=
( t^k + t^{-k} ) - ( t^{k-1} + t^{1-k} ) + \cdots + (-1)^{k},
$
and the result follows.
\EPf
%%%%%%%%%%%%%%%%%%%%%%%%%%%%%%%%%%%%%%%%%%%%%%%%%%%%%%%שש
%\pn
\begin{remark}
This congruence may be used as a simple criterion to prove that some knots
cannot be two-bridge knots.
There is a more efficient criterion by Kanenobu \cite{Kan2,St00} using the Jones and Q polynomials.
\end{remark}
\pn
We also deduce an analogous result for two-component links (see also \cite[p.~186]{Bu})
\begin{corollary}[Modulo 2 Hosokawa polynomials of two-bridge links]\label{hosokawa}
Let $\Delta (t)= \bigl( t^{\undemi} - t^{- \undemi} \bigr) \Bigl(
 a_0 -a_1 ( t+ t^{-1}) + a_2 ( t^2 + t^{-2} ) - \cdots +(-1)^n a_n ( t^n + t^{-n} ) \Bigr)$
be the Alexander polynomial of a two-component two-bridge link. Then all the coefficients $a_i$ are even or
there exists an integer $k\le n$ such that $ a_k, a_{k-2}, a_{k-4}, \ldots$ are odd, and the other coefficients are even.
\end{corollary}
\Pf
If $K$ is a two-component two-bridge link, its Conway polynomial is an odd Fibonacci
polynomial modulo 2, that is of the form $ f_{2h} (z)$.
An easy induction shows that
$$ f_{4k} \bigl( t^{\undemi} - t^{- \undemi} \bigr)=
\bigl( t^{\undemi} - t^{- \undemi} \bigr) \bigl( u_1+u_3+ \cdots + u_{2k-1} \bigr)$$
and
$$ f_{4k+2} \bigl( t^{\undemi} - t^{- \undemi} \bigr)=
\bigl( t^{\undemi} - t^{- \undemi} \bigr) \bigl( 1+ u_2+ \cdots + u_{2k} \bigr),$$
where $u_j= t^j + t^{-j},$ and the result follows.
\EPf
Theorem \ref{alphabis} implies both  Murasugi and  Hartley theorems for two-bridge knots.
\begin{thm}[Murasugi (1958), Hartley (1979)] {\label{al2}}
Let
$$\Delta(t) = a_0 -a_1 (t+ t^{-1}) + a_2 (t^2+t^{-2}) - \cdots +(-1)^n a_n (t^n + t^{-n})
, \  a_n>0 $$
be the Alexander polynomial of a two-bridge knot. There exists an integer $k \le n$ such that
$ a_0=a_1=\ldots =a_k > a_{k+1}> \ldots > a_n$.
\end{thm}
\Pf
Let $K$ be a two-bridge knot and
$\nabla(z)= \alpha_0 f_1 - \alpha _1 f_3 + \cdots + (-1)^n \alpha_{n} f_{2n+1}$
be its Conway polynomial written in the Fibonacci basis.
By Theorem \ref{alphabis},  $\alpha _n \alpha _k \ge 0$ for all $k$, and  if $\alpha _i=0$ for some
$i$ then $\alpha_j=0$ for $j\le i$.

Let
$\Delta (t) = a_0 -a_1 ( t+ t^{-1}) + a_2 ( t^2 + t^{-2} ) - \cdots +(-1)^n a_n ( t^n + t^{-n} )
, \   a_n>0$
be the Alexander polynomial of $K$.
We have $ \Delta (t)= \eps \nabla ( t^{ \undemi} - t^{- \undemi} ) $,  where $\eps= \pm 1,$
and then,
by  Corollary \ref{chvar}, $  \eps \alpha_k = a_k- a_{k+1}$.
We deduce that $ \eps \alpha_n= a_n >0,$ and then
$ a_k- a_{k+1} = \eps \alpha_k \ge 0 $ for all $k$.
Consequently we obtain $ a_0  \ge a_1 \ge \ldots  \ge  a_n >0$.
\pn
Furthermore, if $a_k= a_{k-1}$ for some $k,$ then $\alpha_{k-1}=0,$ and consequently
$ \alpha_{j-1}=0$ for all $j\le k$.
This implies that for all $j \le k, $  $a_j=a_{j-1} $, which concludes the proof.
\EPf
%%%%%%%%%%%%%%%%%%%%%%%%%%%%%%%%ששששששש
\pn
Now, we shall give explicit formulas for  Alexander coefficients in terms of
 Conway coefficients.
\begin{proposition}\label{c2a}
Let $ Q(z)={\tilde c}_{0} + {\tilde c}_{1} z^2 + \cdots + {\tilde c}_{n} z^{2n}$ be a polynomial.
 We have
$$
Q(t^{\undemi} - t^{- \undemi}) =
 a_0 -a_1 ( t+ t^{-1}) + a_2 ( t^2 + t^{-2} ) - \cdots +(-1)^n a_n ( t^n + t^{-n} ),
$$
where
\[
a_{n-j} =  \sum _{k=0}^j (-1)^{n-k}{\tilde c}_{n-k} \comb{2n-2k}{j-k}.\label{c2af}
\]
\end{proposition}
\Pf
It is sufficient to prove Formula (\ref{c2af}) for the monomials
$Q(z)=z^{2m}$.
Let us consider $u_i=t^i + t^{-i}. $ By the binomial formula we have
$$
 \Bigl( t^{\undemi} - t^{- \undemi} \Bigr)^{2m}=
\sum _{k=0} ^{m-1} (-1)^k \comb{2m}{k} u_{m-k} +
(-1)^m \comb{2m}{m}.
$$
and then
$a_{n-j}= (-1)^m \comb {2m}h $ where $ m-h=n-j.$
On the other hand, the proposed formula asserts
$$
a_{n-j} =  \sum _{k=0}^j (-1)^{n-k}{\tilde c}_{n-k} \comb{2n-2k}{j-k} =
(-1)^m \comb {2m}h \quad {\rm where} \quad h=m+j-n,
$$
which is the same result.
\EPf
%%%%%%%%%%%%%%%%%%%%%%%%%%%%%%%%%%%%%%%%%%%%%%%%%%%
\begin{remark}
Considering the Fibonacci polynomials $f_{2n+1} = \sum_{k=0}^n \comb{2n-k}{k} z^{2n-2k}$,
Formulas (\ref{afib}) and (\ref{c2af}) give
the identity
$$
1=\sum_{k=0}^j (-1)^{k} \comb{2n-k}{k} \comb{2n-2k}{j-k}, \, n,j \geq 0.
$$
\end{remark}
\begin{remark}\label{ra2cf}
Fukuhara  \cite{Fu} gives a converse  formula for the $c_k$ in terms of the $a_k$,
$$
{\tilde c}_{n-j} = \sum _{k=0} ^j (-1)^{n-k} a_{n-k} \frac{2n-2k}{2n-j-k} \comb{2n-j-k}{2n-2j}.
%\label{a2cf}
$$
\end{remark}
From the bounds we obtained for  Conway coefficients we
can deduce a simple proof  of the Nakanishi--Suketa bounds (\cite[Th. 1, 2]{NS})
for  the Alexander  coefficients.
\begin{thm}[Nakanishi--Suketa (1993)]\label{naks}
We have the following sharp inequalities
(where all the $a_i$ are positive):
\bn
\item
$a_{n-j} \le a_n \Bigl(  \sum _{k=0} ^j \comb{2n-2k}{j-k}\comb{ 2n-k}{k} \Bigr)$.
\item
$ 2 a_n -1 \le a_{n-1} \le (4n-2) a_n +1$.
\en
\end{thm}
\Pf

\bn
\item
Using Formula (\ref{c2af}) and Theorem \ref{nak1bis}, we obtain
\[
\abs{a_{n-j}} \le
\sum_ {k=0}^j \abs{\tilde c_{n-k} } \comb{2n-2k}{j-k} \le
\abs{a_n} \sum_{k=0}^j \comb{2n-k}{k}\comb{2n-2k}{j-k}.\label{ns1}
\]
\item
We have
$ \abs {\tilde c_{n-1}} \le \comb{2n-2}1 \abs { \tilde c_n}  +1$ by Theorem \ref{nak3bis},
and $a_{n-1} = \tilde c_{n-1} - \comb{2n}1 \tilde c_n$ by Proposition \ref{c2a}. We thus deduce
\[
\abs {a_{n-1}} \le \comb{2n}1 \abs {\tilde c_n} + \comb{2n-2}1 \abs{ \tilde c_n} +1=
(4n-2) \abs{ a_n} +1.\label{ns1b}
\]
We also have
$$
\abs{ a_{n-1}} \ge \comb{2n}1 \abs{\tilde c_{n}} - \abs{ \tilde c_{n-1}} \ge
 \comb{2n}1 \abs{\tilde c_n} - \comb{2n-2}1 \abs{\tilde c_n} -1 = 2 \abs{a_n} -1.
$$
\en
The upper bounds (\ref{ns1}) and (\ref{ns1b}) are attained by the knots $ C(2,2, \ldots, 2)$.
\EPf
%%%%%%%%%%%%%%%%%%%%%%%%%%%%%%%%%%%%%%%%
\pn
We also have the following sharp bound, which improves
the Nakanishi--Suketa third bound (\cite[Th. 3]{NS})
\begin{thm}\label{nakss}
If $a_n \not = 1,$ then
$a_{n-2} \le ( 8n^2-15n+8) a_n + 2n-1$. This bound is sharp.
\end{thm}
\Pf
From Proposition \ref{c2a} and Theorem \ref{nak3bis}, we get
\begin{align*}
\abs{a_{n-2}} &\le
\comb {2n}2 \abs{ \tilde c_n} +
\comb {2n-2} 1 \abs{ \tilde c_{n-1}} +
\comb {2n-4}0 \abs {\tilde c_{n-2} }\\
&\le \comb {2n}2 \abs{ \tilde c_n}
+ \comb{2n-2}1 (\comb{2n-2}1 \abs { \tilde c_n}  +1)
+ \Bigl(\comb{2n-3}{2} + \frac 1g ( \comb{2n-3}{1} -1)\Bigr)\abs { \tilde c_n}  +1\\
&=
(8n^2-16n+10 + \frac{2(n-2)}g    ) \abs{a_n} + 2n-1 .
%\Bigl( 1 + (2n-2)^2 + \comb{2n-3}{2} +\frac 1g ( \comb{2n-3}{1} -1)\Bigr)\abs { \tilde c_n} + 2n-1,
\end{align*}
If $a_n \ne 1$ then $g \ge 2,$ and we obtain
\[
\abs{ a_{n-2}} \le \abs{ a_n} ( 8n^2 -15n +8) + 2n-1. \label{ns2}
\]
This bound  is attained for the knot $C(4,2,2,2, \ldots, 2)$.
\EPf
\begin{example}
Let us consider the Conway polynomial $ \nabla_ K (z)= 1 +8 z^2 +3 z
^4 -z^6$ of the knot $K=13n1862$ (see \cite{KA}).
It does not verify the bounds of theorem \ref{nak1bis}, and then it is not a two-bridge knot.
Nevertheless, its Alexander polynomial
$ \Delta _K (t)= 23 - 19 (t+1/t) +9 (t^2+1/t^2) - (t^3+1/t^3) $
 satisfies the bounds of Nakanishi and Suketa,   and also the conditions
of Murasugi and Hartley.
This example shows that the conditions on the Conway coefficients are stronger than
 the conditions on the Alexander coefficient deduced from them.
\end{example}
\begin{remarks}
\bn
\item If $g\geq 3$, the inequality (\ref{ns2}) can be improved:
$$
a_{n-2} \leq (8n^2-16n+10+ \frac{2(n-2)}g  ) a_n + 2n-1.
$$
\item For $j=3$ we obtain
\begin{align*}
a_{n-3}  &\le
2/3\, \left( 2\,n-3 \right)  \left( 8\,{n}^{2}-24\,n+25 \right) a_n+{
\frac { \left( 3\,n-5 \right)  \left( 2\,n-5 \right)}{g}  a_n}+n \left( 2
\,n-3 \right)\\
&\le
1/6\, \left( 64\,{n}^{3}-270\,{n}^{2}+413\,n-225 \right) a_n+n \left( 2
\,n-3 \right).
\end{align*}
\item Since the inequalities on  Conway coefficients are simpler and stronger,
we shall not give the inequalities on Alexander coefficients for $j \ge 4$.
Furthermore, if we want to apply our bounds to the Alexander polynomials, we first compute
$$
{\tilde c}_{n-j} = \sum _{k=0} ^j (-1)^{n-k} a_{n-k} \frac{2n-2k}{2n-j-k} \comb{2n-j-k}{2n-2j},
$$
using Remark \ref{ra2cf} and test if
$\abs{{\tilde c}_{n-j}} \leq \comb{2n-j}{j}   \abs {\tilde c_n}$, which is stronger than the inequality (\ref{ns1}),
or if $\abs{{\tilde c}_{n-j}} \le \Bigl( \comb{ 2n-j-1}{j}
 + \frac 1g \bigl( \comb{ 2n-j-1}{j-1} -1 \bigr)
\Bigr) \abs {c_n} +1$.
The cost of these evaluations is less than the cost of the evaluations of
the  inequalities of Theorem \ref{naks}. They are also sharper.
\en
\end{remarks}
The following example shows an infinite family of polynomials
satisfying all the necessary conditions except the equality case of Theorem \ref{nak1bis}.
\begin{example}
Consider the polynomial  $P(z)= f_{m+1} (z) -2dz^2, \  m =4n \ge 4, \  d \ne 0$.
All its coefficients, except one, satisfy $ c_{m-2k}= { m -k \choose k}$.
By Theorem \ref{nak1bis}, it is not the Conway polynomial of a two-bridge knot.
Hence, the corresponding Alexander polynomial
$$\Delta(t) = 4d+1 -(2d+1)u_1 + u_2 -u_3 + \cdots +  u_{2n},$$
where
$u_i= t^i + t^{-i}$, is not the Alexander polynomial of a two-bridge knot.
Nevertheless, it satisfies all the necessary conditions of Hartley and Murasugi.
If $ 0< d < \frac 12 n(n+1)$, it also satisfies the bounds of
Theorems \ref{nak1bis} and \ref{nak3bis}, and
then the Nakanishi--Suketa bounds.
\end{example}
\section{Conjectures}
We observed a trapezoidal property for the Conway polynomials of two-bridged links with
20 or fewer crossings (their number is $131\,839$).
\begin{conjecture}\label{trapezoidal}
Let
$\nabla_m = c_m \Bigl(\sum_{i=0}^{\pent m2} (-1)^i \alpha_{i}
f_{m+1-2i}\Bigr), \ \alpha_0=1,$ be the Conway polynomial of a
two-bridge link (or knot) written in the Fibonacci basis.
Then there exists $n \leq \pent m2$ such that
$$0 \leq \alpha_0 \leq \alpha_1 \le \cdots \leq \alpha_n, \quad
 \alpha_n  \geq \alpha_{n+1} \geq \cdots \geq \alpha_{\pent m2} \geq 0.$$
\end{conjecture}
\pn
If this conjecture was true, it would imply the following property of
Alexander polynomials:
\begin{conjecture}\label{con-con}
Let $\Delta(t) = a_0 -a_1 ( t+ t^{-1}) + a_2 ( t^2 + t^{-2} ) - \cdots
+(-1)^n a_n ( t^n + t^{-n} )$ be the Alexander polynomial of a
two-bridge knot. Then there exists an integer $k \le n$ such that
$ (a_0, \ldots, a_k) $ is convex and $(a_{k}, \ldots, a_n)$ is concave.
\end{conjecture}
It is shown in \cite{NS} that the sequence $a_j$ is not convex.
\pn
The following conjecture is attributed to Hoste:
\begin{conjecture}[Hoste]\label{hoste}
If $z \in \mathbf{C}$ is a root of the
Alexander polynomial of an alternating knot,
then $\Re z > -1$.
\end{conjecture}
This conjecture is shown to be true in some peculiar cases (see \cite{LM,St11}).
As a direct consequence of the definition of Euler continuant polynomials, we show that:
\begin{theorem}
Let $K$ be a two-bridge link (or knot). Let $\alpha$ be a root of the Alexander polynomial $\Delta_K$, then
$-\Frac{3}{2}<\Re \alpha <3 + 2\sqrt 2$.
If $\alpha$ is real then $3-2\sqrt 2 < \alpha < 3+2\sqrt 2$.
\end{theorem}
\Pf Let $K$ be a two-bridge link. $\nabla_K$ is an Euler continuant polynomial $D_m(b_1, \ldots, b_m)$. If $z$ is a root of $\nabla_K$, then the determinant in Formula (\ref{tridiagonal}) is equal to 0.
It is a classical result in linear algebra that there exists $i$ such that $\abs{b_i z} < 2$. We thus deduce that $\abs{z}<2$.

Let $\alpha$ be a root of $\Delta_K$. Then  $z = \alpha^{\undemi}  - \alpha^{-\undemi} $ is a root of $\nabla_K$ and %therefore $\abs z < 2$.
we have the relation $ P(\alpha,z) = \alpha^2 - (z^2+2) \alpha +1 =0$.
Eliminating $z$ between $P$ and $\abs{z}<2$, we obtain that
$\alpha=x+iy$ satisfies $R(x,y)<0$ where
$$
R =
{x}^{4}+2\,{x}^{2}{y}^{2}+{y}^{4}-4\,{x}^{3}-4\,x{y}^{2}-10\,{x}^{2}-
14\,{y}^{2}-4\,x+1.
$$
\begin{figure}[th]
\begin{center}
{\scalebox{.3}{\includegraphics{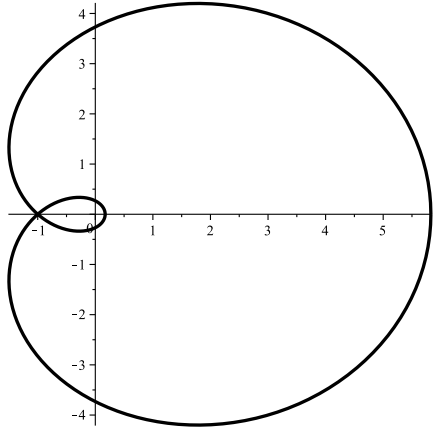}}}
\end{center}
\caption{Region $(R<0)$ containing the roots of Alexander polynomials of two-bridge links.}
\label{alex1}
\end{figure}
An easy computation shows that the curve $R=0$ has vertical tangents at the four points:
$$
(-\Frac 32, \pm \Frac{\sqrt 7}2), \, (3 \pm 2\sqrt 2,0).
$$
Suppose that $\alpha$ is real. Then $z^2=\alpha+1/\alpha-2$ is real and
$\mathrm{Discr}(P) = z^2(z^2+4)\geq 0$.
We thus deduce that $z$ is real and belongs to $(-2,2)$. We thus have
$\alpha \in (3-2\sqrt 2,3+2\sqrt 2)$.
\EPf
\pn
This result is an improvement of those obtained in \cite{LM}.
We found that it was independently obtained by Stoimenow (see \cite{St13}).
It should be improved by a careful study of the tridiagonal matrix $A_m$ in Formula (\ref{tridiagonal}).
\bibliographystyle{line}

\vfill
\pn
\hrule width 7cm height 1pt %depth 3pt
\pn
{\small
Pierre-Vincent {\sc Koseleff}\\
Universit\'e Pierre et Marie Curie (UPMC Sorbonne Univesit\'es),\\
Institut de Math\'ematiques de Jussieu (IMJ-PRG)  \& Inria-Rocquencourt\\
%4, place Jussieu, F-75252 Paris Cedex 05 \\
e-mail: {\tt koseleff@math.jussieu.fr}
\pn
Daniel {\sc Pecker}\\
Universit\'e Pierre et Marie Curie (UPMC Sorbonne Univesit\'es),\\
Institut de Math\'ematiques de Jussieu (IMJ-PRG),\\
%4, place Jussieu, F-75252 Paris Cedex 05 \\
e-mail: {\tt pecker@math.jussieu.fr}
%\noindent MSC2000: 57M25
%\label{LastPage}
\end{document}